\documentclass[10pt,twoside]{article}
\usepackage{amsmath,amsthm,latexsym,amsfonts,amssymb}
\usepackage{Latex-document}
\title{\bf The Topology of {\boldmath $Out(F_n)$}\vskip 6mm}
  \author{Mladen Bestvina\thanks{Department of
Mathematics, University of Utah, USA. E-mail:
bestvina@math.utah.edu}\vspace*{-0.5cm}}
\date{\vspace{-8mm}}
\newtheorem{thm}{Theorem}

\newtheorem{cor}[thm]{Corollary}

\theoremstyle{remark}
\newtheorem{defn}[thm]{\bf Definition}

\def\R{{\mathbb R}}
\def\Z{{\mathbb Z}}

\def\vcd{\operatorname{vcd}}

\begin{document}
\maketitle

\thispagestyle{first} \setcounter{page}{373}

\begin{abstract}

\vskip 3mm

We will survey the work on the topology of $Out(F_n)$ in the last 20 years or so. Much of the development is
driven by the tantalizing analogy with mapping class groups. Unfortunately, $Out(F_n)$ is more complicated and
less well-behaved.

Culler and Vogtmann constructed Outer Space $X_n$, the analog of Teichm\" uller space, a contractible complex on
which $Out(F_n)$ acts with finite stabilizers. Paths in $X_n$ can be generated using ``foldings'' of graphs, an
operation introduced by Stallings to give alternative solutions for many algorithmic questions about free groups.
The most conceptual proof of the contractibility of $X_n$ involves folding.

There is a normal form of an automorphism, analogous to Thurston's normal form for surface homeomorphisms. This
normal form, called a ``(relative) train track map'', consists of a cellular map on a graph and has good
properties with respect to iteration. One may think of building an automorphism in stages, adding to the previous
stages a building block that either grows exponentially or polynomially. A complicating feature is that these
blocks are not ``disjoint'' as in Thurston's theory, but interact as upper stages can map over the lower stages.

Applications include the study of growth rates (a surprising feature of free group automorphisms is that the
growth rate of $f$ is generally different from the growth rate of $f^{-1}$), of the fixed subgroup of a given
automorphism, and the proof of the Tits alternative for $Out(F_n)$. For the latter, in addition to train track
methods, one needs to consider an appropriate version of ``attracting laminations'' to understand the dynamics of
exponentially growing automorphisms and run the ``ping-pong'' argument. The Tits alternative is thus reduced to
groups consisting of polynomially growing automorphisms, and this is handled by the analog of Kolchin's theorem
(this is one instance where $Out(F_n)$ resembles $GL_n(\Z)$ more than a mapping class group).
\\
\indent Morse theory has made its appearance in the subject in several guises. The original proof of the
contractibility of $X_n$ used a kind of ``combinatorial'' Morse function (adding contractible subcomplexes one at
a time and studying the intersections). Hatcher-Vogtmann developed a ``Cerf theory'' for graphs. This is a
parametrized version of Morse theory and it allows them to prove homological stability results. One can
``bordify'' Outer Space (by analogy with the Borel-Serre construction for arithmetic groups) to make the action of
$Out(F_n)$ cocompact and then use Morse theory (with values in a certain ordered set) to study the connectivity at
infinity of this new space. The result is that $Out(F_n)$ is a virtual duality group.
\\
\indent Culler-Morgan have compactified Outer Space, in analogy with Thurston's compactification of Teichm\" uller
space. Ideal points are represented by actions of $F_n$ on $\R$-trees. The work of Rips on group actions on
$\R$-trees can be used to analyze individual points and the dynamics of the action of $Out(F_n)$ on the boundary.
The topological dimension of the compactified Outer Space and of the boundary have been computed. The orbits in
the boundary are not dense; however, there is a unique minimal closed invariant set. Automorphisms with
irreducible powers act on compactified Outer Space with the standard North Pole -- South Pole dynamics. By first
finding fixed points in the boundary of Outer Space, one constructs a ``hierarchical decomposition'' of the
underlying free group, analogous to the Thurston decomposition of a surface homeomorphism.
\\
\indent The geometry of Outer Space is not well understood. The most promising metric is not even symmetric, but
this seems to be forced by the nature of $Out(F_n)$. Understanding the geometry would most likely allow one to
prove rigidity results for $Out(F_n)$.
\\
\vskip 4.5mm

\end{abstract}

\vskip 12mm

\section{Introduction} \setzero

\vskip-5mm \hspace{5mm}

The aim of this note is to survey some of the topological methods
developed in the last 20 years to study the group $Out(F_n)$ of outer
automorphisms of a free group $F_n$ of rank $n$. For an excellent and
more detailed survey see also \cite{MRkaren}. Stallings' paper
\cite{MR813103} marks the turning point and for the earlier history of
the subject the reader is referred to \cite{MR58:28182}. $Out(F_n)$ is
defined as the quotient of the group $Aut(F_n)$ of all automorphisms
of $F_n$ by the subgroup of inner automorphisms. On one hand,
abelianizing $F_n$ produces an epimorphism $Out(F_n)\to
Out(\Z^n)=GL_n(\Z)$, and on the other hand $Out(F_n)$ contains as a
subgroup the mapping class group of any compact surface with
fundamental group $F_n$. A {\it leitmotiv} in the subject, promoted by
Karen Vogtmann, is that $Out(F_n)$ satisfies a mix of properties, some
inherited from mapping class groups, and others from arithmetic
groups. The table below summarizes the parallels between topological
objects associated with these groups.

\begin{table}[th]
\begin{center}
\begin{tabular}{|l|l|l|l|}\hline
         {\footnotesize\bf Mapping}
        &{\footnotesize\bf $Out(F_n)$}
        &{\footnotesize\bf $GL_n(\Z)$}
        &{\footnotesize\bf algebraic }
\\
         {\footnotesize\bf class groups}
        &
        &{\footnotesize\bf (arithmetic groups)}
        &{\footnotesize\bf properties}
\\\hline\hline
         {\footnotesize\sf Teichm\"uller}
        &{\footnotesize\sf Culler-Vogtmann's}
        &{\footnotesize\sf $GL_n(\R)/O_n$}
        &{\footnotesize\sf finiteness properties}
\\
         {\footnotesize\sf space}
        &{\footnotesize\sf Outer space}
        &{\footnotesize\sf (symmetric spaces)}
        &{\footnotesize\sf cohomological dimension}
\\\hline
         {\footnotesize\sf Thurston}
        &{\footnotesize\sf train track}
        &{\footnotesize\sf Jordan}
        &{\footnotesize\sf growth rates}
\\
         {\footnotesize\sf normal form}
        &{\footnotesize\sf representative}
        &{\footnotesize\sf normal form}
        &{\footnotesize\sf fixed points (subgroups)}
\\\hline
         {\footnotesize\sf Harer's}
        &{\footnotesize\sf bordification of}
        &{\footnotesize\sf Borel-Serre}
        &{\footnotesize\sf Bieri-Eckmann}
\\
         {\footnotesize\sf bordification}
        &{\footnotesize\sf Outer space}
        &{\footnotesize\sf bordification}
        &{\footnotesize\sf duality}
\\\hline
         {\footnotesize\sf measured}
        &{\footnotesize\sf {$\R$}-trees}
        &{\footnotesize\sf flag manifold}
        &{\footnotesize\sf Kolchin theorem}
\\
         {\footnotesize\sf laminations}
        &
        &{\footnotesize\sf (Furstenberg boundary)}
        &{\footnotesize\sf Tits alternative}
\\\hline
         {\footnotesize\sf Harvey's}
        &{\footnotesize\sf ?}
        &{\footnotesize\sf Tits}
        &{\footnotesize\sf rigidity}
\\
         {\footnotesize\sf curve complex}
        &
        &{\footnotesize\sf building}
        &
\\\hline
\end{tabular}
\end{center}
\end{table}
Outer space is not a manifold and only a polyhedron, imposing a combinatorial character on $Out(F_n)$.

\section{Stallings' Folds} \setzero
\vskip-5mm \hspace{5mm}

A {\it graph} is a 1-dimensional cell complex. A map $f:G\to G'$
between graphs is {\it simplicial} if it maps vertices to vertices and
open 1-cells homeomorphically to open 1-cells. The simplicial map $f$
is a {\it fold} if it is surjective and identifies two edges that
share at least one vertex. A fold is a homotopy equivalence unless the
two edges share both pairs of endpoints and in that case the induced
homomorphism in $\pi_1$ corresponds to killing a basis element.

\begin{thm}[Stallings \cite{MR85m:05037a}]\label{stallings}
A simplicial map $f:G\to G'$ between finite connected graphs can be
factored as the composition $$G=G_0\to G_1\to \cdots\to G_k\to G'$$
where each $G_i\to G_{i+1}$ is a fold and $G_k\to G'$ is locally
injective (an immersion). Moreover, such a factorization can be found
by a (fast) algorithm.
\end{thm}

In the absence of valence 1 vertices the last map $G_k\to G'$ can be
thought of as the core of the covering space of $G'$ corresponding to
the image in $\pi_1$ of $f$. The following problems can be solved
algorithmically using Theorem \ref{stallings} (these were known
earlier, but Theorem \ref{stallings} provides a simple unified argument). Let
$F$ be a free group with a fixed finite basis.
\vskip -1mm
\begin{itemize}
\item Find a basis of the subgroup $H$ generated by a given finite
collection $h_1,\cdots,h_k$ of elements of $F$.
\item Given $w\in F$, decide if $w\in <h_1,\cdots,h_k>$.
\item Given $w\in F$, decide if $w$ is conjugate into
$<h_1,\cdots,h_k>$.
\item Given a homomorphism $\phi:F\to F'$ between two free groups of
finite rank, decide if $\phi$ is injective, surjective.
\item Given finitely generated $H<F$ decide if it has finite index.
\item Given two f.g. subgroups $H_1,H_2<F$ compute $H_1\cap H_2$ and
also the collection of subgroups $H_1\cap H_2^g$ where $g\in F$. In
particular, is $H_1$ malnormal?
\item Represent a given automorphism of $F$ as the composition of
generators of $Aut(F)$ of the following form:
\vskip 0.1cm
\subitem Signed permutations: each $a_i$ maps to $a_i$ or to $a_i^{-1}$.
\vskip 0.1cm
\subitem Change of maximal tree: $a_1\mapsto a_1$, $a_i\mapsto
a_1^{\pm 1}a_i$ or $a_i\mapsto a_ia_1^{\pm 1}$ ($i>1$).
\item Todd-Coxeter process \cite{MR88g:20084}.
\end{itemize}

\section{Culler-Vogtmann's Outer space} \setzero
\vskip-5mm \hspace{5mm}

Fix the wedge of $n$ circles $R_n$ and a natural identification
$\pi_1(R_n)\cong F_n$ in which oriented edges correspond to the basis
elements. Thus any $\phi\in Out(F_n)$ can be thought of as a homotopy
equivalence $R_n\to R_n$. A {\it marked metric graph} is a pair
$(G,g)$ where
\begin{itemize}
\item $G$ is a finite graph without vertices of valence 1 or 2.
\item $g:R_n\to G$ is a homotopy equivalence (the {\it marking}).
\item $G$ is equipped with a path metric so that the sum of the
lengths of all edges is 1.
\end{itemize}

{\it Outer space} $X_n$ is the set of equivalence classes of marked
metric graphs under the equivalence relation $(G,g)\sim (G',g')$ if
there is an isometry $h:G\to G'$ such that $gh$ and $g'$ are
homotopic \cite{MR87f:20048}.

If $\alpha$ is a loop in $R_n$ we have the length function
$l_\alpha:X_n\to \R$ where $l_\alpha(G,g)$ is the length of the
immersed loop homotopic to $g(\alpha)$. The collection $\{l_\alpha\}$
as $\alpha$ ranges over all immersed loops in $R_n$ defines an
injection $X_n\to\R^\infty$ and the topology on $X_n$ is defined so
that this injection is an embedding. $X_n$ naturally decomposes into
open simplices obtained by varying edge-lengths on a fixed marked
graph. The group $Out(F_n)$ acts on $X_n$ on the right via
$$(G,g)\phi=(G,g\phi).$$

\begin{thm}[Culler-Vogtmann \cite{MR87f:20048}]
$X_n$ is contractible and the action of $Out(F_n)$ is properly
discontinuous (with finite point stabilizers). $X_n$ equivariantly
deformation retracts to a $(2n-3)$-dimensional complex ($n>1$).
\end{thm}

If $(G,g)$ and $(G',g')$ represent two points of $X_n$, there is a
``difference of markings'' map $h:G\to G'$ such that $hg$ and $g'$ are
homotopic. Representing $h$ as a composition of folds (appropriately
interpreted) leads to a path in $X_n$ from $(G,g)$ to
$(G',g')$. Arranging that these paths vary continuously with endpoints
leads to a proof of contractibility of $X_n$
\cite{MRsteiner},\cite{MRskora},\cite{MR93i:57002}.

\begin{cor} The virtual cohomological dimension
$\vcd(Out(F_n))=2n-3$ ($n>1$).\end{cor}

\begin{thm}[Culler \cite{MR86g:20027}]\label{culler}
Every finite subgroup of $Out(F_n)$ fixes a point of $X_n$.
\end{thm}

Outer space can be equivariantly compactified
\cite{MR88f:20055}. Points at infinity are represented by actions of
$F_n$ on $\R$-trees.

\section{Train tracks}
\vskip-5mm \hspace{5mm}

Any $\phi\in Out(F_n)$ can be represented as a cellular map $f:G\to G$
on a marked graph $G$. We say that $\phi$ is {\it reducible} if there
is such a representative where
\begin{itemize}
\item $G$ has no vertices of valence 1 or 2, and
\item there is a proper $f$-invariant subgraph of $G$ with at least
one non-contractible component.
\end{itemize}
Otherwise, we say that $\phi$ is {\it irreducible}.

A cellular map $f:G\to G$ is a {\it train track map} if for every
$k>0$ the map $f^k:G\to G$ is locally injective on every open
1-cell. For example, homeomorphisms are train track maps and Culler's
theorem guarantees that every $\phi\in Out(F_n)$ of finite order has a
representative $f:G\to G$ which is a homeomorphism. More generally, we
have

\begin{thm}[Bestvina-Handel \cite{MR92m:20017}]\label{bh}
Every irreducible outer automorphism $\phi$ can be represented as a
train track map $f:G\to G$.
\end{thm}

Any vertex $v\in G$ has a cone neighborhood, and the frontier points
can be thought of as ``germs of directions'' at $v$. A train track map (or
any cellular map that does not collapse edges) $f$ induces the
``derivative'' map $Df$ on these germs (on possibly different
vertices). We declare two germs at the same vertex to be equivalent
(and the corresponding ``turn'' {\it illegal}) if they get identified
by some power of $Df$ (and otherwise the turn is {\it legal}). An immersed
loop in $G$ is {\it legal} if every turn determined by
entering and then exiting a vertex is legal. It follows that
$f$ sends legal loops to legal loops. This gives a method for
computing the growth rate of $\phi$, as follows. The {\it transition
matrix} $(a_{ij})$ of $f$ (or more generally of a cellular map $G\to G$ that is
locally injective on edges) has $a_{ij}$ equal to the number of times
that the $f$-image of $j^{th}$ edge crosses $i^{th}$ edge.
Applying the
Perron-Frobenius theorem to the transition matrix, one can find a
unique metric structure on $G$ such that $f$ expands lengths of edges
(and also legal loops) by a factor $\lambda\geq 1$. For a conjugacy
class $\gamma$ in $F_n$ the growth rate is defined as
$$GR(\phi,\gamma)=\limsup_{k\to\infty}\log(||\phi^k(\gamma)||)/k$$
where $||\gamma||$ is the word length of the cyclically reduced word
representing $\gamma$. Growth rates can be computed using lengths of
loops in
$G$ rather than in $R_n$.
\begin{cor}
If $\phi$ is irreducible as above, then either $\gamma$ is a
$\phi$-periodic conjugacy class, or
$GR(\phi,\gamma)=\log\lambda$. Moreover, $\limsup$ can be replaced by
$\lim$.
\end{cor}

The proof of Theorem \ref{bh} uses a folding process that successively
reduces the Perron-Frobenius number of the transition matrix until
either a train track representative is found, or else a reduction of
$\phi$ is discovered. This process is algorithmic (see
\cite{MR96d:57014},\cite{MR2001e:57017}).

Another application of train tracks is to fixed subgroups.

\begin{thm}[Bestvina-Handel \cite{MR92m:20017}]\label{scott}
Let $\Phi:F_n\to F_n$ be an automorphism whose associated outer
automorphism is irreducible. Then the fixed subgroup $Fix(\Phi)$ is
trivial or cyclic. Without the irreducibility assumption, the rank of
$Fix(\Phi)$ is at most $n$.
\end{thm}

It was known earlier by the work of Gersten \cite{MR88f:20042} that
$Fix(\Phi)$ has finite rank (for simpler proofs see
\cite{MR87m:20096},\cite{MR89a:20024}).  The last sentence in the
above theorem was conjectured by Peter Scott. Subsequent work by
Collins-Turner \cite{MR96j:20037}, Dicks-Ventura \cite{MR94i:20047},
Ventura \cite{MR98h:20042}, Martino-Ventura \cite{MR2001h:20029},
imposed further restrictions on a subgroup of $F_n$ that occurs as the
fixed subgroup of an automorphism. To analyze reducible automorphisms,
a more general version of a train track map is required.

\begin{defn}\it
A cellular map $f:G\to G$ on a finite graph with no vertices of
valence 1 that does not collapse any edges is a {\it relative train
track map} if there is a filtration
$$\emptyset=G_0\subset\cdots\subset G_m=G$$ into $f$-invariant
subgraphs with the following properties. Denote by $H_r$ the closure
of $G_r\setminus G_{r-1}$, and by $M_r$ the part of the transition
matrix corresponding to $H_r$. Then $M_r$ is the zero matrix or an
irreducible matrix. If $M_r$ is irreducible and the Perron-Frobenius
eigenvalue $\lambda_r>1$ then:
\begin{itemize}
\item the derivative $Df$ maps the germs in $H_r$ to germs in $H_r$,
\item if $\alpha$ is a nontrivial path in $G_{r-1}$ with endpoints in
$G_{r-1}\cap H_r$ then $f(\alpha)$, after pulling tight, is also a
nontrivial path with endpoints in $G_{r-1}\cap H_r$, and
\item every legal path in $H_r$ is mapped to a path that does not
cross illegal turns in $H_r$.
\end{itemize}
\end{defn}\rm

As an example, consider the automorphism $a\mapsto a, b\mapsto ab,
c\mapsto caba^{-1}b^{-1}d$, $d\mapsto dbcd$ represented on the rose
$R_4$. The strata are $\emptyset\subset G_1=\{a\}\subset \{a,b\}\subset
G$. $H_1$ and $H_2$ have $\lambda=1$ while $H_3$ has $\lambda_3>1$.
The following is an analog of Thurston's normal form for surface
homeomorphisms.

\begin{thm}\label{bh2}{\rm\cite{MR92m:20017}}
Every automorphism of $F_n$ admits a relative train track
representative.
\end{thm}

Consequently, automorphisms of $F_n$ can be thought of as being built
from building blocks (exponential and non-exponential kinds) but the
later stages are allowed to map over the previous stages. This makes
the study of automorphisms of $F_n$ more difficult (and interesting)
than the study of surface homeomorphisms. Other non-surface
phenomena (present in linear groups) are:
\begin{itemize}
\item stacking up non-exponential strata produces
(nonlinear) polynomial growth,
\item the growth rate of an automorphism is generally different from
the growth rate of its inverse.
\end{itemize}

\section{Related spaces and structures}
\vskip-5mm \hspace{5mm}

Unfortunately, relative train track representatives are far from
unique. As a replacement, one looks for canonical objects associated
to automorphisms that can be computed using relative train
tracks. There are 3 kinds of such objects, all stemming from the
surface theory: laminations, $\R$-trees, and hierarchical
decompositions of $F_n$ \cite{MR97f:20047}.

\noindent
{\bf Laminations.}  Laminations were used in the proof
of the Tits alternative for $Out(F_n)$. To each automorphism one
associates finitely many attracting laminations. Each consists of a
collection of ``leaves'', i.e. biinfinite paths in the graph $G$, or
alternatively, of an $F_n$-orbit of pairs of distinct points in the
Cantor set of ends of $F_n$. A leaf $\ell$ can be computed by
iterating an edge in an exponentially growing stratum $H_r$. The other
leaves are biinfinite paths whose finite subpaths appear as subpaths
of $\ell$. Some of the attracting laminations may be sublaminations of
other attracting laminations, and one focuses on the maximal (or {\it
topmost}) laminations. It is possible to identify the basin of
attraction for each such lamination. Let $\cal H$ be any subgroup of
$Out(F_n)$.  Some of the time it is possible to find to elements
$f,g\in\cal H$ that attract each other's laminations and then the
standard ping-pong argument shows that $<f,g>\cong
F_2$.  Otherwise, there is a finite set of attracting laminations
permuted by $\cal H$, a finite index subgroup ${\cal H}_0\subset \cal
H$ that fixes each of these laminations and a homomorphism (``stretch
factor'') ${\cal H}_0\to A$ to a finitely generated abelian group $A$
whose kernel consists entirely of polynomially growing
automorphisms. There is an analog of Kolchin's theorem that says that
finitely generated groups of polynomially growing automorphisms can
simultaneously be realized as relative train track maps on the same
graph (the classical Kolchin theorem says that a group of unipotent
matrices can be conjugated to be upper triangular, or equivalently
that it fixes a point in the flag manifold). The main step in the
proof of the analog of Kolchin's theorem is to find an appropriate
fixed $\R$-tree in the boundary of Outer space. This leads to the Tits
alternative for $Out(F_n)$:

\begin{thm} [Bestvina-Feighn-Handel
\cite{MR98c:20045},\cite{MR2002a:20034},\cite{MRbfh:book2}]
Any subgroup $\cal H$ of $Out(F_n)$ either contains $F_2$ or is
virtually solvable.
\end{thm}

A companion theorem \cite{MRbfh:book3} (for a simpler proof see
\cite{MRemina}) is that solvable subgroups of $Out(F_n)$ are virtually
abelian.

\noindent
{\bf $\R$-trees.}
Points in the compactified Outer space are represented as
$F_n$-actions on $\R$-trees. It is then not surprising that the Rips
machine \cite{MR96h:20056}, which is used to understand individual
actions, provides a new tool to be deployed to study $Out(F_n)$.
Gaboriau, Levitt, and Lustig \cite{MR99d:20035} and Sela
\cite{MR97f:20047} find another proof of Theorem \ref{scott}. Gaboriau
and Levitt compute the topological dimension of the boundary of Outer
Space \cite{MR97c:20039}. Levitt and Lustig show \cite{MRll} that
automorphisms with irreducible powers have the standard north-south
dynamics on the compactified Outer space. Guirardel
\cite{MR2002f:20061} shows that the action of $Out(F_n)$ on the
boundary does not have dense orbits; however, there is a unique
minimal closed invariant set. For other applications of $\R$-trees in
geometric group theory, the reader is referred to the survey
\cite{MRsurvey}.

\noindent
{\bf Cerf theory.}
An advantage of $Aut(F_n)$ over $Out(F_n)$ is that there is a natural
inclusion $Aut(F_n)\to Aut(F_{n+1})$. One can define {\it Auter Space}
$AX_n$ similarly to Outer space, except that all graphs are equipped
with a base vertex, which is allowed to have valence 2. The degree of
the base vertex $v$ is $2n-\mbox{valence}(v)$. Denote by $D^k_n$ the
subcomplex of $AX_n$ consisting of graphs of degree $\leq
k$. Hatcher-Vogtmann \cite{MR2000e:20041} develop a version of Cerf
theory and show that $D^k_n$ is $(k-1)$-connected. Since the quotient
$D_n^k/Aut(F_n)$ stabilizes when $n$ is large, one sees that
(rational) homology $H_i(Aut(F_n))$ also stabilizes when $n$ is large
($n\geq 3i/2$). Hatcher-Vogtmann show that the same is true for
integral homology and in the range $n\geq 2i+3$. They also make
explicit computations in low dimensions \cite{MR99m:20127} and
all stable rational homology groups $H_i$ vanish for $i\leq 7$.

\noindent
{\bf Bordification.}
The action of $Out(F_n)$ on Outer space $X_n$ is not cocompact. By
analogy with Borel-Serre bordification of symmetric spaces
\cite{MR52:8337} and Harer's bordification of Teichm\" uller space
\cite{MR87c:32030}, Bestvina and Feighn \cite{MR2001m:20041} bordify
$X_n$, i.e. equivariantly add ideal points so that the action on the
new space $BX_n$ is cocompact. This is done by separately
compactifying every simplex with missing faces in $X_n$ and then
gluing these together. To see the idea, consider the case of the
theta-graph in rank 2. Varying metrics yields a 2-simplex $\sigma$
without the vertices. As a sequence of metrics approaches a missing
vertex, the lengths of two edges converge to 0. Restricting a metric
to these two edges and normalizing so that the total length is 1 gives
a point in $[0,1]$ (the length of one of the edges), and a way to
compactify $\sigma$ by adding an interval for each missing vertex. The
compactified $\sigma$ is a hexagon. This procedure equips the limiting
theta graph with a metric that may vanish on two edges, in which case
a ``secondary metric'' is defined on their union. In general, a graph
representing a point in the bordification is equipped with a sequence
of metrics, each defined on the core of the subgraph where the
previous metric vanishes.

Lengths of curves (at various scales) provide a ``Morse function'' on
$BX_n$ with values in a product of $[0,\infty)$'s with the target
lexicographically ordered. The sublevel and superlevel sets intersect
each cell in a semi-algebraic set and it is possible to study how the
homotopy types change as the level changes. A distinct advantage of
$BX_n$ over the spine of $X_n$ (an equivariant deformation retract) is
that the change in homotopy type of superlevel sets as the level
decreases is very simple -- via attaching of cells of a fixed dimension.

\begin{thm}[Bestvina-Feighn \cite{MR2001m:20041}]
$BX_n$ and $Out(F_n)$ are $(2n-5)$-connected at infinity, and
$Out(F_n)$ is a virtual duality group of dimension $2n-3$.
\end{thm}

\noindent
{\bf Mapping tori.}
If $\phi:F_n\to F_n$ is an automorphism, form the mapping torus
$M(\phi)$. This is the fundamental group of the mapping torus $G\times
[0,1]/(x,1)\sim (f(x),0)$ of any representative $f:G\to G$, and it
plays the role analogous to 3-manifolds that fiber over the circle. Such a
group is always coherent \cite{MR2000i:20050}. A quasi-isometry
classification of these groups seems out of reach, but the following
is known. When $\phi$ has no periodic conjugacy classes, $M(\phi)$ is
a hyperbolic group \cite{MR2001m:20061}. When $\phi$ has polynomial
growth, $M(\phi)$ satisfies quadratic isoperimetric inequality
\cite{MR2001k:20089} and moreover, $M(\phi)$ quasi-isometric to
$M(\psi)$ for $\psi$ growing polynomially forces $\psi$ to grow as a
polynomial of the same degree \cite{MRnatasa}. Bridson and Groves
announced \cite{MRbg} that $M(\phi)$ satisfies quadratic isoperimetric
inequality for all $\phi$.

\noindent
{\bf Geometry.}  Perhaps the biggest challenge in the field is to find
a good geometry that goes with $Out(F_n)$. The payoff would most
likely include rigidity theorems for $Out(F_n)$. Both mapping class
groups and arithmetic groups act isometrically on spaces of
nonpositive curvature. Unfortunately, the results to date for
$Out(F_n)$ are negative. Bridson \cite{MR94c:57040} showed that Outer
space does not admit an equivariant piecewise Euclidean $CAT(0)$
metric. $Out(F_n)$ ($n>2$) is far from being $CAT(0)$
\cite{MR96j:20050},\cite{MR94j:20043}.

An example of a likely rigidity theorem is that higher rank
lattices in simple Lie groups do not embed into $Out(F_n)$. A possible
strategy is to follow the proof in \cite{MRbefu} of the analogous fact
for mapping class groups. The major missing piece of the puzzle is the
replacement for Harvey's curve complex; a possible candidate is
described in \cite{MR99i:20038}.



\providecommand{\bysame}{\leavevmode\hbox to3em{\hrulefill}\thinspace}
\providecommand{\MR}{\relax\ifhmode\unskip\space\fi MR }
\providecommand{\MRhref}[2]{%
  \href{http://www.ams.org/mathscinet-getitem?mr=#1}{#2}
}
\providecommand{\href}[2]{#2}

\label{lastpage}

\end{document}